\documentclass[12pt,reqno]{amsart}
\newtheorem{theorem}{\bf Theorem}

\newtheorem{lemma}[theorem] {\bf Lemma}

\newcommand{\NN}{\mathbb N}

\newcommand{\RR}{\mathbb R}
\newcommand{\CC}{\mathbb C}

\DeclareMathOperator{\re}{Re\,\!}
\DeclareMathOperator{\im}{Im\,\!}

\usepackage{geometry,amssymb}
\begin{document}
\title {On the real linear polarization constant problem}

\author[M. Matolcsi and G. A. Mu{\~n}oz]
{M{\'a}t{\'e} Matolcsi and Gustavo A. Mu{\~n}oz}

\address[M{\'a}t{\'e} Matolcsi]{A. R{\'e}nyi Institute for Mathematics, Hungarian Academy of Sciences, Budapest, P.O.B. 127, 1364 Hungary.}
\email{matomate@renyi.hu}

\address[Gustavo A. Mu{\~n}oz]{Facultad de Matem{\'a}ticas, Departamento de An{\'a}lisis \\
Universidad Complutense de Madrid, Madrid 28040, Spain.}
\email{gustavo\_fernandez@mat.ucm.es}

\thanks{The first author was supported by OTKA-T047276.}
\thanks{ The second author was supported by BFM 2003-06420.}

\begin{abstract}
The present paper deals with lower bounds for the norm of products
of linear forms. It has been proved by J. Arias-de-Reyna
\cite{ARIAS}, that %for ${\mathbb C}^n$,
the so-called $n^{\rm th}$ linear polarization constant
$c_n({\mathbb C}^n)$ is $n^{n/2}$, for arbitrary $n\in\NN$. The
same value for $c_n({\mathbb R}^n)$ is only conjectured. In a
recent work A. Pappas and S. R{\'e}v{\'e}sz prove that $c_n({\mathbb
R}^n)=n^{n/2}$ for $n \le 5$. Moreover, they show that if the
linear forms are given as $f_j(x)=\langle x,a_j \rangle$, for some
unit vectors $a_j$ $(1\leq j\leq n)$, then the product of the
$f_j$'s attains at least the value $n^{-n/2}$ at the normalized
signed sum of the vectors $\{a_j\}_{j=1}^{n}$ having maximal length. Thus
they asked whether this phenomenon remains true for arbitrary
$n\in{\mathbb N}$. We show that for vector systems
$\{a_j\}_{j=1}^{n}$ close to an orthonormal system, the
Pappas-R{\'e}v{\'e}sz estimate does hold true. Furthermore, among these
vector systems the only system giving $n^{-n/2}$ as the norm of
the product is the orthonormal system. On the other hand, for
arbitrary vector systems we answer the question of A. Pappas and
S. R{\'e}v{\'e}sz in the negative when $n\in {\mathbb N}$ is large enough.
We also discuss various further examples and counterexamples that
may be instructive for further research towards the determination
of $c_n(\RR^n)$.
\end{abstract}

\maketitle

\section{\bf Introduction and notation}

For convenience we recall the basic definitions needed to discuss
polynomials on a Banach space. If ${\mathbb K}$ is the real or
complex field and $X$ is a Banach space over ${\mathbb K}$, then
we denote by $B_X$ and $S_X$ the closed unit ball and the unit
sphere of $X$ respectively. A map $P:X\rightarrow {\mathbb K}$ is
a (continuous) {\em $n$-homogeneous polynomial} if there is a
(continuous) symmetric $n$-linear mapping $L:X^n \rightarrow
{\mathbb K}$ for which $P(x)=L(x,\ldots,x)$ for all $x\in X$. In
this case it is convenient to write $P={\widehat L}$. We let
${\mathcal P}(^{n}X)$ denote the space of scalar valued continuous
$n$-homogeneous polynomials on $X$. We define the norm of a
(continuous) homogeneous polynomial $P:X\rightarrow {\mathbb K}$
by
    $$
    \|P\|=\sup\{|P(x)|:\ x\in B_X\}.
    $$
For general background on polynomials, we refer to \cite{DINEEN}.

If $P_k\in {\mathcal P}(^{n_k}X)$ ($1\leq k\leq m$) then the
pointwise product of the $P_k$'s given by $(P_1\cdots
P_m)(x):=P_1(x)\cdots P_m(x)$ for every $x\in X$ is also a
homogeneous polynomial, in fact if $n=n_1+\ldots+n_m$ we have
$P_1\cdot P_2\cdots P_m\in {\mathcal P}(^nX)$. Clearly
$\|P_1\cdots P_m\|\leq \|P_1\|\cdots \|P_m\|$. An estimate in the
other direction is more difficult to establish, but still
possible. Indeed, there is an absolute constant
$C_{n_1,\ldots,n_m}$ such that
    \begin{equation}\label{prodpol}
    \|P_1\|\cdot\|P_2\|\cdots\|P_m\|\leq
    C_{n_1,\ldots,n_m}\|P_1\cdot P_2\cdots P_m\|.
    \end{equation}
Products of polynomials and estimates like \eqref{prodpol} have
been studied by several authors. For a general account on this
problem we recommend \cite{BST} and the references therein.

In this paper we will restrict ourselves to the case where the
$P_k$'s are continuous linear functionals $f_1, f_2, \ldots, f_n$
on $X$. Then the product $(f_1f_2 \cdots f_n)(x):= f_1(x)f_2(x)
\cdots f_n(x)$ is a continuous $n$-homogeneous polynomial on $X$
and from (\ref{prodpol}) there exists $C_n>0$ such that
    \begin{equation}\label{prodlinforms}
    \|f_1\|\|f_2\| \cdots \|f_n\| \leq C_n \|f_1f_2 \cdots f_n \|\,.
    \end{equation}
Estimate (\ref{prodlinforms}) was also studied by R.A. Ryan and B.
Turett in \cite{RT}. In \cite{BST} it was proved, in the case of
{\it complex} Banach spaces, that $C_n \leq n^n$ and the constant
$n^n$ is best possible in general. However $n^n$ can be improved
for specific spaces. The best fitting constant in
\eqref{prodlinforms} was defined by C. Ben{\'\i}tez, Y. Sarantopulos
and A. Tonge in \cite{BST} as
    $$
    c_n(X):= \inf \{M >0: \|f_{1}\| \cdots \|f_{n}\| \leq M \cdot
    \|f_{1} \cdots f_{n}\|, \forall f_{1}, \ldots ,f_{n} \in
    X^{\ast}\},
    $$
and it is referred to in the literature as the $n^{th}$ linear
polarization constant of $X$.

Let us represent the Hilbert space of the $n$-tuples of elements
of ${\mathbb K}$ endowed with the Euclidean norm $\|\cdot\|_2$ by
${\mathbb K}^n$. Then, it has been proved by S. R{\'e}v{\'e}sz and Y.
Sarantopoulos \cite{REVSAR} that
    $$
    c_n(X) \geq c_n({\mathbb K}^n),\quad\forall {n \in {\mathbb
    N}},
    $$
for any infinite dimensional Banach space $X$. This inequality
shows the importance of $c_n({\mathbb K}^n)$ when estimating
$c_n(X)$, at least for infinite dimensional Banach spaces.

J. Arias-de-Reyna proved in \cite{ARIAS} that $c_n({\mathbb
C}^n)=n^{\frac{n}{2}}$, however his proof does not adapt to the
real case. It has been conjectured in \cite{BST} that
$c_n({\mathbb R}^n)=n^{\frac{n}{2}}$ also holds, but no proof has
been given yet.

Observe that in order to determine $c_n({\mathbb R}^n)$ one only
needs to consider norm-one functionals $f_1,f_2,\ldots,f_n$ in
(\ref{prodlinforms}). Therefore by the Riesz Representation
Theorem a polynomial $P_{n}(x):=f_1(x)\cdots f_n(x)$, where $f_k
\in S_{({{\mathbb R}^n})^{\ast}}$, $1 \leq k \leq n$, can be
written in the form
\begin{equation}\label{Pndef}
P_{n}(x)=\langle x,a_{1} \rangle \cdot \langle x,a_{2} \rangle
\cdot \ldots \langle x,a_{n} \rangle\;,
\end{equation}
where $a_{j} \in S_{({{\mathbb R}^n})^{\ast}}=S_{{{\mathbb
R}^n}}$.

If ${\mathcal B}_n=\{e_j:\ 1\leq j\leq n\}$ is the canonical basis
of ${\mathbb R}^n$ and we put $a_j=e_j$ for $1\leq j\leq n$ in
(\ref{Pndef}), then for $x\in B_{{\mathbb R}^n}$ with coordinates
$(x_1,\ldots,x_n)$, we have (using the fact that the geometric
mean is smaller than the quadratic mean):
    \begin{align*}
    |P_n(x)|=|x_1\cdots x_n|\leq
    n^{-\frac{n}{2}}(x_1^2+\ldots+x_n^2)^{\frac{n}{2}}=n^{-\frac{n}{2}}\|x\|_2^n\leq n^{-\frac{n}{2}},
    \end{align*}
from which $\|P_n\|\leq n^{-\frac{n}{2}}$. In addition to this
    $$
    |P_n(\frac{1}{\sqrt{n}},\ldots,\frac{1}{\sqrt{n}})|=n^{-\frac{n}{2}},
    $$
and therefore $\|P_n\|=n^{-\frac{n}{2}}$. This lets us state that
$c_n({\mathbb R}^n)\geq n^{\frac{n}{2}}$. In order to prove the
reverse inequality one could try and find for every $P_n$ as in
(\ref{Pndef}) a norm one vector $x$ in ${\mathbb R}^n$ satisfying
$|P_n(x)|\geq n^{-\frac{n}{2}}$. This has been done in the real
case for $n\leq 5$ by A. Pappas and S. R{\'e}v{\'e}sz in \cite{PR}, taking
$x$ as the normalization of the signed combination of the
functional vectors $\{a_j\}_{j=1}^n$ with maximal length (see next
section).

A. Pappas and S. R{\'e}v{\'e}sz also asked whether their argument can be
generalized for any dimension. We show in Section 4 that their
argument fails in spaces of large enough dimension.

Another approach to the problem of estimating $c_n({\mathbb R}^n)$
is using complexification arguments together with the result of J.
Arias-de-Reyna in $\CC^n$. In fact, this idea has been used by S.
R{\'e}v{\'e}sz and Y. Sarantopoulos in \cite{REVSAR} to prove
    \begin{equation}\label{bestestimate}
    c_n({\mathbb R}^n)\leq 2^{\frac{n}{2}-1}n^{\frac{n}{2}},
    \end{equation}
which is the best known estimate on $c_n({\mathbb R}^n)$ (see also
\cite{MST} for a complete account on results on complexifications
of polynomials).

With the notation $\|P_n\|_{S_{{\mathbb K}^n}}:=\sup\{|P_n(\xi)|:\
\xi\in S_{{\mathbb K}^n}\}$ one could aim to prove that
$\|P_n\|_{S_{{\mathbb R}^n}}=\|P_n\|_{S_{{\mathbb C}^n}}$ for any
$\{a_j\}_{j=1}^n\subset S_{{\mathbb R}^n}$. Then it would follow
that
    $$
    \|P_n\|_{S_{{\mathbb R}^n}}=\|P_n\|_{S_{{\mathbb C}^n}}\leq n^{-\frac{n}{2}},
    $$
from which $c_n({\mathbb R}^n)=n^{\frac{n}{2}}$. However, it is
shown in Section \ref{secComp} that $\|P_n\|_{S_{{\mathbb R}^n}}$
does not necessarily coincide with $\|P_n\|_{S_{{\mathbb C}^n}}$.
Moreover, we prove that using complexification arguments it is not
possible to improve the estimate \eqref{bestestimate}.

\section{\bf Mean vectors of maximal length}

In the following we will refer to a real choice of signs
$\epsilon$ as an $n$-tuple,
$\epsilon=(\epsilon_1,\ldots,\epsilon_n)$, with $\epsilon_j=\pm 1$
$(1\leq j\leq n)$. If $a_1,\dots,a_n\in S_{\RR^n}$ are $n$ vectors
in ${\mathbb R}^n$ we define
    \begin{equation}\label{signedcomb}
    a_\epsilon:=\sum_{i=1}^n \epsilon_i a_i.
    \end{equation}
If we select $\epsilon$ to maximize $\|a_\epsilon\|_2$ it can be
easily shown that $\left< a_\epsilon , \epsilon_j a_j \right> \ge
1$ for $1\leq j\leq n$. Indeed, if we fix $j$ ($1\leq j\leq n$)
and $\epsilon^\prime$ is the choice of signs given by
$\epsilon_k^\prime=\epsilon_k$ if $k\ne j$ and
$\epsilon_j^\prime=-\epsilon_j$ then
    $$
    \|a_\epsilon\|_2^2\geq \|a_{\epsilon^\prime}\|_2^2=
    \|a_{\epsilon}\|^2+4(1-\left< a_\epsilon , \epsilon_j a_j\right>),
    $$
from which $\left< a_\epsilon , \epsilon_j a_j \right> \ge 1$
follows immediately. Note that if we replace the $a_j$'s in
(\ref{Pndef}) by $\epsilon_ja_j$, the norm of $P_n$ does not
change. Therefore, we can assume without loss of generality that
the choice $\epsilon =(1,1,\dots ,1)$ gives the maximal length,
and then, by the argument above we have
    \begin{align}
    y_1: & = \langle a_{1},a_{1} \rangle+ \langle a_{1},a_{2} \rangle+
    \ldots + \langle a_{1},a_{n} \rangle \geq 1\nonumber\\
    y_2: & =\langle a_{2},a_{1} \rangle+ \langle a_{2},a_{2} \rangle+
    \ldots + \langle a_{2},a_{n} \rangle \geq 1\nonumber\\
      & ..................................................\label{coordinates}\\
    y_n: & =\langle a_{n},a_{1} \rangle+ \langle a_{n},a_{2} \rangle+
    \ldots + \langle a_{n},a_{n} \rangle \geq 1.\nonumber
    \end{align}
A. Pappas and S. R{\'e}v{\'e}sz used in \cite{PR} this type of signed
combinations. In particular they considered the normalized mean
vector
\begin{equation}\label{amean}
x:=\frac {a_\epsilon}{\|a_\epsilon\|_2}= \frac
{a_1+\cdots+a_n}{\|a_1+\cdots+a_n\|_2},
\end{equation}
(with $a_\epsilon$ having maximal length), and they obtained the
following:
\begin{theorem}[A. Pappas and S. R{\'e}v{\'e}sz \cite{PR}]\label{lowcases}
Let $n\le 5$. If $P_n$ is as in (\ref{Pndef}) and $x$ as in
(\ref{amean}) then
\begin{equation}\label{qu}
|P_n(x)|=\prod_{j=1}^n |\langle x\,,\, a_j\rangle| \ge n^{-n/2}\,.
\end{equation}
\end{theorem}

\vspace{0.2in}

\noindent {\bf Remarks.} (1) Let us note that \eqref{coordinates}
is the special case $r_j=1$ $(j=1,\dots,n)$ of a more general
result, known as Bang's Lemma, see \cite{BANG}. For the special
case of $r_j=1$ above, the argument occurs at several places, see
eg. \cite[Lemma 2.4.1 (i)]{LMS}, \cite{MM} and \cite{PR}.

(2) We observe that signed combinations of vectors such as
(\ref{signedcomb}) have been used in many other constructions (see
for instance \cite{GV} and \cite{LMS}). In all cases the usual
approach is to choose the combination of maximal length.

\vspace{0.2in}

\section{\bf A Question of A. Pappas and S. R{\'e}v{\'e}sz and its complex version}

Estimate (\ref{qu}) entails the conjectured value $n^{n/2}$ of the
polarization constant for dimensions $n\le 5$. That is why the
next question was posed by A. Pappas and S. R{\'e}v{\'e}sz \cite{PR}.

\vskip 0.10in

\noindent {\bf Question} (A. Pappas and S. R{\'e}v{\'e}sz \cite{PR}). Is
it true, that for any $n$ and unit vectors $a_j\in\RR^n$ $(1\leq
j\leq n)$, with signs of the unit vectors $a_j$ $(1\leq j\leq n)$
chosen to maximize the length of the signed sum (and thus
satisfying \eqref{coordinates}, too), the mean value vector
\eqref{amean} satisfies $|P(x)|\ge n^{-n/2}$?

\vskip 0.10in

This question has a natural analogue in the complex case. Instead
of $\pm 1$ we consider complex choices of sign, which we define by
$n$-tuples $\epsilon=(\epsilon_1,\ldots,\epsilon_n)$, where
$\epsilon_j=e^{i\varphi_j}$ for some $\varphi_j\in [0,2\pi)$
($1\leq j\leq n$). In this setting we consider complex vectors
$a_j\in S_{\CC^n}$ $(j=1,\dots,n)$ and in this case $a_\epsilon$
is defined as
    $$
    a_\epsilon:=\sum_{j=1}^n\epsilon_ja_j=\sum_{j=1}^ne^{i\varphi_j}a_j.
    $$
If we choose $e^{i\varphi_j}$ ($1\leq j\leq n$) so as to maximize
the norm of $a_\epsilon$, as in the real case it can be easily
shown that
    $$
    \langle a_\epsilon,\epsilon_ja_j\rangle=
    \langle\sum_{k=1}^ne^{i\varphi_k}a_j,e^{i\varphi_j}a_j\rangle\geq
    1,
    $$
for $1\leq j\leq n$. Indeed, taking $1\leq j\leq n$,
$u:=\sum_{k\ne j} e^{i\varphi_k}a_k$, $v:=e^{i\varphi_j}a_j$ and
$\lambda\in {\mathbb C}$ such that $|\lambda|=1$ and $\langle
u,\lambda v\rangle=\left|\langle u, v\rangle\right|$, we have
    \begin{align*}
    \|u+v\|_2^2\geq \|u+\lambda v\|_2^2 & \Leftrightarrow
    \|u\|_2^2+\|v\|_2^2+2\re (\langle u,v\rangle)\geq
    \|u\|_2^2+\|v\|_2^2+2\re (\langle u,\lambda v\rangle)\\
    & \Leftrightarrow \re (\langle u,v\rangle)\geq \re (\langle u,\lambda
    v\rangle)= |\langle u,v\rangle|.
    \end{align*}
Therefore $\im(\langle u,v\rangle)=0$ and $\langle
u,v\rangle=|\langle u,v\rangle|\geq 0$. This implies the desired
result:
    $$
    \langle a_\epsilon,\epsilon_ja_j\rangle=
    \langle\sum_{k=1}^ne^{i\varphi_k}a_j,e^{i\varphi_j}a_j\rangle=
    \langle u+v,v\rangle=\langle u,v\rangle+1\geq 1.
    $$
Now, since
    $$
    |P_n(z)|=\left|\langle a_1,z\rangle \dots \langle a_n,z\rangle
    \right|=\left|\langle e^{i\varphi_1}a_1,z\rangle \dots \langle
    e^{i\varphi_n}a_n,z\rangle\right|,
    $$
we can replace the $a_j$'s in (\ref{Pndef}) by $e^{i\varphi_j}a_j$
without loss of generality. In other words, writing $y_j:= \langle
\sum_{k=1}^n a_k , a_j\rangle$, we find $y_j\ge 1$, $j=1,\dots,n$,
as in (\ref{coordinates}). From here, the argument proving Theorem
\ref{lowcases} coincides with the real case (cf. \cite{PR}),
considering the normalized mean vector $z$:
    $$
    z:=\frac{a_1+\ldots+a_n}{\|a_1+\ldots+a_n\|_2}.
    $$
The arguments of \cite{PR} (applied verbatim) then show that $z$
always satisfies $|P_n(z)| \ge n^{-n/2}$, for $n\le 5$.

The advantage of this result (in these low dimensional cases)
compared to \cite{ARIAS} and \cite{BALL1} is that it gives a {\it
construction} for the vector $z$ as opposed to proving existence
only.

The analogue to the question of A. Pappas and S. R{\'e}v{\'e}sz in the
complex case is to ask whether this construction of $z$ works in
{\it all} dimensions.

\section{\bf Mean vectors of maximal length for systems
close to orthonormal vectors}

It is plausible to expect that the only system of vectors for
which $\|P_n\|\le n^{-n/2}$ is the orthonormal system (which
satisfies $\|P_n\|= n^{-n/2}$ as shown in Section 1). Even in the
complex case, this statement does not follow from the
considerations in \cite{ARIAS} and \cite{BALL1}, and remains an
open question. In this section we prove the {\em local} uniqueness
of the orthonormal system in both the real and the complex case.
In Theorem \ref{closeon} below we show that for vectors close to
the orthonormal system we always have $\|P_n\|\ge n^{-n/2}$ and
the Pappas-R{\'e}v{\'e}sz choice \eqref{amean} of mean vector belonging to
maximal length always provides strict inequality (unless the
system is orthonormal, in which case strict inequality is not
possible). With the help of some specific examples, however, we
will show later that the Pappas-R{\'e}v{\'e}sz choice \eqref{amean} of
mean vector $x$ does not satisfy $|P_n(x)|\ge n^{-n/2}$ in
general.

\begin{theorem}\label{closeon} Let $H$ be any real or complex Hilbert
space and $\{a_j\}_{j=1}^n\subset S_H$. We assume (without loss of
generality) that the choice of signs $\epsilon =(1,1, \dots ,1)$
gives maximal length among the mean vectros $a_{\epsilon}$. (This
also means that the system $\{a_j\}_{j=1}^n$ satisfies condition
\eqref{coordinates}.) Assume also, with the notation used in
\eqref{coordinates}, that $y_j\le 3.5$, for $j=1,\dots,n$. Then
$\|P_n\|\ge n^{-n/2}$ holds with equality only when $y_j=1$ for
$j=1,\dots,n$, that is, only for the orthonormal vector system.
\end{theorem}

\begin{proof}
Before proceeding with the proof we note that the condition
$y_j\le 3.5$ is certainly satisfied in a small neighbourhood of
the orthonormal system. (Due to the dimension being finite, all
definitions of 'neighbourhood' are equivalent. Vector systems in a
small neighbourhood are obtained by small perturbations of the
vectors of an orthonormal system.)

Let $x$ be the mean vector defined as in \eqref{amean}. The
assertion $|P_n(x)|\ge n^{-n/2}$ is equivalent to state that the
inequality
    \begin{equation}\label{yineq}
    y_1^2y_2^2\cdots y_n^2 \ge \left(\frac{y_1+y_2+\cdots
    +y_n}{n}\right)^n
    \end{equation}
holds true (see \cite{PR}). With a slight reformulation, we put
$t_i:=y_i-1$ $(i=1,\dots,n)$ and write
    \begin{equation}\label{tineq}
    (1+t_1)^2(1+t_2)^2\cdots (1+t_n)^2 \ge
    \left(1+\frac{t_1+t_2+\cdots +t_n}{n}\right)^n\,.
    \end{equation}
Note that here $0\le t_i\le n-1$ for all $i=1,\dots,n$. Now it
suffices to show \eqref{tineq} for real quantities $0\le t_j \le
2.5$ ($j=1,\dots,n$). To start with, observe that
    \begin{equation}\label{expineq}
    (1+t)^2\ge e^{t}.
    \end{equation}
is satisfied for $0\le t \le t_0$, where $t_0$ is the (unique)
root of $f(t):=(1+t)^2-e^t$ in $(1,\infty)$. Moreover, since
$t_0>2.5$, \eqref{expineq} follows with strict inequality if
$0<t\le 2.5$.

Multiplying together this inequality for all $t_i$ with
$i=1,\dots,n$, we obtain
    \begin{equation}\label{texpineq}
    (1+t_1)^2(1+t_2)^2\cdots (1+t_n)^2 \ge e^{t_1+t_2+\cdots+t_n}\,.
    \end{equation}
Since $e^x\ge (1+x/n)^n$ for all $x\ge 0$ and $n\in \NN$, we
conclude \eqref{tineq}, and hence the assertion.

We have equality in \eqref{expineq} only for $t=0$, hence for
equality in \eqref{texpineq} we must have $t_1=\dots=t_n=0$. That
is, $\{a_j\}_{j=1}^n$ must be an orthonormal system of vectors.
\end{proof}

The proof of the theorem above heavily relied on the assumption
$y_j\le 3.5$. A thorough analysis of the proof above leads us to
give a negative answer to the question of A. Pappas and S. R{\'e}v{\'e}sz
in high dimensions.

\begin{theorem}\label{nogoodmean}
If $n$ is large enough, then there exist vectors
$\{a_j\}_{j=1}^n\subset S_{\RR^n}$ so that taking the mean vector
\eqref{amean} of maximal length, we have
    \begin{equation}\label{qunot}
    |P_n(x)|=\prod_{j=1}^n |\langle x\,,\, a_j\rangle| < n^{-n/2}\,.
    \end{equation}
The analogous results holds for $a_j\in S_{\CC^n}$ $(j=1,\dots,n)$
and complex signs $e^{i\varphi_j}$.
\end{theorem}

\begin{proof}
The proof is a specific example obtained by analyzing the proof of
Theorem \ref{closeon} above.

Take $n\geq 34$ and let $a_1,\ldots,a_n$ be the $n$ unit vectors
in $H$ defined as follows:
    \begin{align*}
    a_j &
    :=b:=(\frac{1}{\sqrt{6}},\stackrel{(6)}{\ldots},\frac{1}{\sqrt{6}},0\ldots,0)
    \quad \textrm{for $1\leq j\leq 6$}\\ a_j & :=e_j \quad\textrm{for
    $6<j \le n$},
    \end{align*}
where $\{e_j:\ 1\leq j\leq n\}$ is the canonical basis of
${\mathbb K}^n$.

Then it is obvious that a choice of signs maximizing the length of
the mean vector (both in the real and complex case) is
$z=6b+\sum_{j=7}^{n}e_j$, and for this vector we have
    $$
    |P_n({z}/{\|z\|})|=\frac{6^6}{(n+30)^\frac{n}{2}}<\frac{1}{n^\frac{n}{2}}\quad
    \textrm{for $n\geq 34$}.
    $$
\end{proof}

\vspace{0.2in}

\noindent {\bf Remarks.} (1) In the example above there exists a
combination of signs such that the corresponding mean vector $x$
satisfies $|P_n(x)|\ge n^{-n/2}$. But this choice of signs {\it
does not correspond} to the sum vector of maximal length. Indeed,
the choice
$x=a_1+a_2+a_3+a_4-a_5-a_6+\sum_{j=7}^{n}a_j=2b+\sum_{j=7}^{n}e_j$
is such a choice.\newline (2) The example used in the previous
proof can be easily generalized by considering natural numbers
$n>d$ and the system of vectors $\{a_j\}_{j=1}^n$ given by
    \begin{align*}
    a_j &
    :=b:=(\frac{1}{\sqrt{d}},\stackrel{(d)}{\ldots},\frac{1}{\sqrt{d}},0\ldots,0)
    \quad \textrm{for $1\leq j\leq d$}\\ a_j & :=e_j \quad\textrm{for
    $d<j \le n$}.
    \end{align*}
In this case $z=db+\sum_{j=d+1}^{n}e_j$ and
    $$
    |P({z}/{\|z\|})|=\frac{d^d}{(n+d^2-d)^\frac{n}{2}}.
    $$
The reader can check easily that for every $d\in {\mathbb N}$ we
can find $n_0\in {\mathbb N}$ with
    $$
    \frac{d^d}{(n+d^2-d)^\frac{n}{2}}<\frac{1}{n^\frac{n}{2}}\quad
    \textrm{for $n\geq n_0$}.
    $$
It is also a simple exercise to see that the smallest possible
value for $n_0$ is $34$ and it is achieved when $d=6$.

\vspace{0.2in}

\section{\bf Norms over $\RR^n$ and over $\CC^n$ }\label{secComp}

Complexifications provide another natural approach towards the
determination of real polarization constants. In this context, it
is natural to ask whether the norm of a real polynomial of the
form \eqref{Pndef} remains the same when considered over the
complex unit ball $S_{\CC^n}$ instead of $S_{\RR^n}$. In view of
the result of Arias-de-Reyna \cite{ARIAS}, this would imply the
conjectured value $c_n({\mathbb R}^n)=n^{n/2}$. However, in this
section we show that if $n\ge 3$ then the norm of the complex
polynomial can be strictly larger than that of the real one.

The main result in this section is based on the following
well-known estimate by G. A. Mu{\~n}oz, Y. Sarantopoulos and A. Tonge
\cite{MST}:
\begin{theorem}\label{complexificationconstant}
If  $P:{\mathbb R}^n\rightarrow {\mathbb R}$ $(n\geq 2)$ is an
$n$-homogeneous polynomial then
    \begin{equation}\label{realvscomplex}
    \|P\|_{S_{\CC^n}} \leq 2^{\frac{n-2}{2}}\| P\|_{S_{\RR^n}}.
    \end{equation}
Moveover, the constant $2^{\frac{n-2}{2}}$ is sharp and equality
in \eqref{realvscomplex} is achieved for the polynomial
$R_n:{\mathbb R}^n\rightarrow {\mathbb R}$ defined by
    $$
    R_n(x_1,x_2,\ldots,x_n):=\re (x_1+ix_2)^n.
    $$
\end{theorem}

Regarding the real linear polarization problem, one could hope to
achieve a better estimate than \eqref{realvscomplex} by
restricting attention to products of functionals instead of {\it
all} $n$-homogeneous polynomials. However, we shall prove in this
section that the polynomials $R_n$ $(n\in {\mathbb N})$ can be
written as a product of functionals as in \eqref{Pndef}. This fact
leads us to conclude that the inequality \eqref{realvscomplex}
cannot be improved even if we consider only polynomials as
described in \eqref{Pndef}.

\vspace{0.2in}

\noindent {\bf Remark.} The polynomials $R_n$ $(n\geq 2)$ where
originally defined in \cite{MST} on ${\mathbb R}^2$. For
simplicity we shall consider the restriction of $R_n$ to ${\mathbb
R}^2$, in other words
    \begin{equation}\label{extremalpolrealvscomplexR2}
    R_n(x,y):=\re (x+iy)^n.
    \end{equation}

\vspace{0.2in}

The following elementary result shall be required in order to
prove that $R_n$ $(n\in {\mathbb N})$ can be factored as the
product of $n$ linear forms:

\begin{lemma}\label{factorizationlemma}
If $P\in {\mathcal P}(^n{\mathbb R}^2)$ satisfies $P(x_0,y_0)=0$
for some $(x_0,y_0)\in {\mathbb R}^2\backslash \{(0,0)\}$, then
there exists $Q\in {\mathcal P}(^{n-1}{\mathbb R}^2)$ such that
    $$
    P(x,y)=(-y_0x+x_0y)Q(x,y),
    $$
for every $(x,y)\in {\mathbb R}^2$.
\end{lemma}

\begin{proof}
Suppose that $x_0\ne 0$ and define $p(t):=P(1,t)$ for all $t\in
{\mathbb R}$. Then $p$ is a real polynomial of degree at most $n$
such that
    $$
    p(\frac{y_0}{x_0})=P(1,\frac{y_0}{x_0})=\frac{1}{x_0^n}P(x_0,y_0)=0.
    $$
Therefore $p(t)=(t-\frac{y_0}{x_0})q(t)$ for some real polynomial
$q(t):=a_{n-1}t^{n-1}+\ldots+a_1 t+a_0$ of degree at most $n-1$.
Hence, if we take $(x,y)\in {\mathbb R}^2$ with $x\ne 0$ we have
    $$
    P(x,y)=x^np(\frac{y}{x})=(\frac{y}{x}-\frac{y_0}{x_0})
    x^nq(\frac{y}{x})=(-y_0x+x_0y)\frac{x^{n-1}}{x_0}q(\frac{y}{x}).
    $$
On the other hand the mapping $Q:{\mathbb R}^2\rightarrow {\mathbb
R}$ given by
    $$
    Q(x,y)=\left\{ \begin{array}{ll}
    \frac{x^{n-1}}{x_0}q(\frac{y}{x}) & \textrm{if $x\ne 0$}\\
    \frac{a_{n-1}}{x_0}y^{n-1} & \textrm{if $x=0$}  \ \ ,
    \end{array}\right.
    $$
is clearly an $(n-1)$-homogeneous polynomial. By continuity it
follows immediately that $P(x,y)=(-y_0x+x_0y)Q(x,y)$ for all
$(x,y)\in {\mathbb R}^2$.
\end{proof}

\begin{theorem}\label{factorizationRn}
If $R_n$ $(n\in {\mathbb N})$ is as in
\eqref{extremalpolrealvscomplexR2} then there exist
$\{a_j\}\subset S_{{\mathbb R}^2}$ and a constant $K\ne 0$ so that
    $$
    R_n(v)=K\langle v,a_1 \rangle\cdots \langle v,a_n \rangle,
    $$
for all $v\in {\mathbb R}^2$.
\end{theorem}

\begin{proof}
Assume first that $n$ is odd and define (identifying ${\mathbb
R}^2$ with ${\mathbb C}$)
    $$
    (x_j,y_j):=\bigg(\sin\frac{j\pi}{n},-\cos\frac{j\pi}{n}\bigg)=e^{i(\frac{j\pi}{n}-\frac{\pi}{2})}\quad\textrm{for
    $j=0,1,\ldots,n-1$}.
    $$
Then, since $2j-n$ is always an odd integer for all odd $n\in
{\mathbb N}$ and $0\leq j\leq n-1$, it follows that
    $$
    R_n(x_j,y_j):=\re e^{in(\frac{j\pi}{n}-\frac{\pi}{2})}=\re
    e^{i(2j-n)\frac{\pi}{2}}=0,
    $$
and hence by Lemma \ref{factorizationlemma} for every $v\in
{\mathbb R}^2$
    $$
    R_n(v)=K\langle v,a_1 \rangle\cdots \langle v,a_n \rangle,
    $$
where $a_j=(\cos\frac{j\pi}{n},\sin\frac{j\pi}{n})$ for $0\leq
j\leq n-1$ and
    $$
    K=\frac{1}{\prod_{j=0}^{n-1}\cos\frac{j\pi}{n}}=(-1)^{\frac{n-1}{2}}2^{n-1}.
    $$
Now suppose that $n$ is even and define
    $$
    (x_j,y_j):=\bigg(\sin\frac{(2j+1)\pi}{2n},-\cos\frac{(2j+1)\pi}{2n}\bigg)=e^{i(\frac{(2j+1)\pi}{2n}-\frac{\pi}{2})}\quad\textrm{for
    $j=0,1,\ldots,n-1$}.
    $$
Then since $2j+1-n$ is always an odd integer for all even $n\in
{\mathbb N}$ and $1\leq j\leq n-1$ we have
    $$
    R_n(x_j,y_j)=\re e^{in(\frac{(2j+1)\pi}{2n}-\frac{\pi}{2})}=\re
    e^{i(2j+1-n)\frac{\pi}{2}}=0.
    $$
and hence by Lemma \ref{factorizationlemma} for every $v\in
{\mathbb R}^2$
    $$
    R_n(v)=K\langle v,a_1 \rangle\cdots \langle v,a_n \rangle,
    $$
where $a_j=(\cos\frac{(2j+1)\pi}{2n},\sin\frac{(2j+1)\pi}{2n})$
for $0\leq j\leq n-1$ and
    $$
    K=\frac{1}{\prod_{j=0}^{n-1}\cos\frac{(2j+1)\pi}{2n}}=(-1)^{\frac{n}{2}}2^{n-1}.
    $$
\end{proof}

\vspace{0.2in}

\noindent {\bf Remark.} If $a_j=(a_j^1,a_j^2)$ $(1\leq j\leq n)$
are the unit vectors obtained in Theorem \ref{factorizationRn} and
we regard the $a_j$'s as vectors in ${\mathbb R}^n$ $(n\geq 2)$ by
setting $a_j:=(a_j^1,a_j^2,0\ldots,0)\in {\mathbb R}^n$ $(1\leq
j\leq n)$, then the $n$-homogeneous polynomial $P_n:{\mathbb
R}^n\rightarrow {\mathbb R}$ defined by
    $$
    P_n(x):=\langle x,a_1\rangle\cdots\langle x,a_n\rangle,
    $$
is of the type \eqref{Pndef} and satisfies
    $$
    \|P_n\|_{S_{{\mathbb C}^n}}=2^{\frac{n-2}{2}}\|P_n\|_{S_{{\mathbb
    R}^n}}.
    $$
This shows that the complexification argument used to prove
estimate \eqref{bestestimate} cannot be improved.

\vspace{0.2in}

\noindent {\bf Acknowledgements:} The authors are grateful to
Professor Szilard R{\'e}v{\'e}sz for his help in the realization of this
paper.

\end{document}